\documentclass[12pt]{article}
\newtheorem{thm}{Theorem}
\newtheorem{prop}[thm]{Proposition}  
\newtheorem{defn}{Definition} \newtheorem{lemma}[thm]{Lemma} 
\newcommand{\vu}{\underline{u}} \newcommand{\vv}{\underline{v}} 
\newcommand{\va}{\underline{a}} \newcommand{\vb}{\underline{b}} 
\newcommand{\vz}{\underline{z}} \newcommand{\ve}{\underline{e}} 
\newcommand{\vo}{\underline{0}} \newcommand{\vt}{\underline{t}} 
\input{diagrams} 

\title{A geometric theory of harmonic and semi-conformal maps}
\author{Anders Kock}
\begin{document}
\maketitle \small Abstract.  We describe for any Riemannian manifold
$M$ a certain scheme $M_L$, lying in between the first and second
neighbourhood of the diagonal of $M$. Semi-conformal maps between Riemannian
manifolds are then analyzed as those maps that preserve $M_L$;
harmonic maps are analyzed as those that preserve the
(Levi-Civita-) mirror image formation inside $M_L$.

\normalsize
\section*{Introduction} For any Riemannian manifold $M$, we describe a
subscheme $M_L \subseteq M\times M$, which encodes information about
as well conformal as harmonic maps out of $M$ in a succinct geometric
way. Thus, a submersion $\phi :M\to N$ between Riemannian manifolds is
semi-conformal (=horizontally conformal) iff $\phi \times \phi$ maps
$M_L$ into $N_L$ (Theorem \ref{semi-c}); and a map $\phi :M\to N$ is
a harmonic map if it ``commutes with mirror image formation for
$M_L$'', where mirror image formation is one of the manifestations
of the Levi-Civita parallelism (derived from the Riemannian metric).
The mirror image preservation property
is best expressed in the set theoretic language for schemes, which we
elaborate on in Section 1. Then it just becomes the statement: for
$(x,z)\in M_L \subseteq M\times M$,
$\phi (z') =(\phi ( z))'$, where the primes denote mirror image
formation in $x$ (respectively in $\phi (x)$). In particular, when
the codomain is ${\bf R}$ (the real line with standard metric), this
characterization of harmonicity reads
$$\phi(z')= 2\phi (x)-\phi (z),$$
that is, $\phi (x)$ equals  the {\em average value} of $\phi (z)$
and $\phi (z')$, for any $z$ with $(x,z)\in M_L$.

The last section deals with {\em harmonic morphisms} between
Riemannian manifolds, meaning harmonic maps which are at the same
time semi-con\-formal.

This paper has some overlap with \cite{Laplace}, but provides a
simplification of the construction of $M_L$, and hence also of the
proofs. Theorems \ref{semi-c} and \ref{three} below are new. A novelty in
the presentation is a systematic use of the $\log$-$\exp$ bijections
that relate the infinitesimal neighbourhoods like $M_L$ with their linearized
version in the tangent bundle.

The first section is partly expository; it tries to present a (rather 
primitive) version of the category of (affine) schemes, and the ``synthetic'' 
language in which we talk about them.

The paper grew out of a talk presented at 
  the 5th conference ``Geometry and Topology of Manifolds'', 
Krynica 2003; I want to thank the organizers for the invitation.

\section{The language of schemes}
  Let $M$ be a smooth manifold. In the ring $C^{\infty}(M\times M)$,
we have the ideal $I$ of functions vanishing on the diagonal
$M\subseteq M\times M$. K\"{a}hler observed that differential 1-forms
on $M$ may be encoded as elements in $I/I^2$ (the module of
K\"{a}hler differentials) (here, $I^2$ is the ideal of functions
vanishing to the second order on the diagonal). Similarly, elements
of $I^2/I^3$ encode quadratic differential forms on $M$. Using the language of
schemes will allow us to discuss elements of $I/I^2$ or of $I^2/I^3$
in a more geometric way. We summarize here what we need about
schemes. First, note that every smooth manifold $M$ gives rise (in a
contravariant way) to a commutative ${\bf R}$-algebra, the ring
$C^{\infty}(M)$ of (smooth ${\bf R}$-valued) functions on it.
Grothendieck's bold step was to think of {\em any} commutative ${\bf
R}$-algebra as the ring of smooth functions on some ``virtual'' {\em
geometric} object $\overline{A}$, the {\em affine scheme} defined by
$A$.  So $A= C^{\infty}(\overline{A})$, by definition, and the
category of affine schemes $Sch$ is just the opposite (dual) of the
category $Alg$ of (commutative ${\bf R}$-)algebras, $$Sch =
(Alg)^{op}.$$ The category of affine schemes contains the category of
smooth manifolds as a full subcategory: to the manifold $M$, associate
the scheme $\overline{C^{\infty}(M)}$ (which we shall not notationally
distinguish from $M$, except for the manifold ${\bf R}$, where we
write $R$ for $\overline{C^{\infty}({\bf R})}$).

Some  important schemes associated to a manifold $M$ are its
infinitesimal
``neighbourhoods of the diagonal'' $M_{(k)}$, considered classically  by
Grothen\-dieck \cite{gr}, Malgrange \cite{ma}, Kumpera and Spencer
\cite{ks} and others. For each natural number $k$, $M_{(k)} \subseteq
M\times M$
is the subscheme of $M\times M$ given by the algebra
$C^{\infty}(M\times M)/I^{k+1}$, where $I$ is the ideal of functions
vanishing on the ``diagonal'' $M\subseteq M\times M$; thus $I^{k+1}$ is 
the ideal of functions vanishing to the $k+1$'st order on the diagonal.

We have $M\subseteq M_{(1)} \subseteq M_{(2)} \subseteq \ldots
\subseteq M\times M$, with $M\subseteq M\times M$ identified with the
submanifold consisting of ``diagonal'' points $(x,x)$.

Now, by definition,
$$C^{\infty}(M\times M)/I^3 = C^{\infty}(M_{(2)}),$$
so in the language of schemes, we arrive at the following way of
speaking: elements of $C^{\infty}(M\times M)/I^3 $ are {\em functions}
on $M_{(2)}$; and elements in $I^2/I^3 \subseteq  C^{\infty}(M\times
M)/I^3 $ are functions on $M_{(2)}$ which {\em vanish} on $M_{(1)}$.

(A similar geometric language was presented in \cite{Krynica} for the
elements of $I/I^2$ (=the K\"{a}hler differentials): they are
functions on $M_{(1)}$ vanishing om $M_{(0)} = M$, i.e.\  they are
{\em combinatorial differential 1-forms} in the sense of \cite{SDG}.)

Synthetic differential geometry adds one feature to this aspect of
scheme theory, namely extended use of set theoretic language for
speaking about objects in (sufficiently nice) categories, like $Sch$.
Thus, since $M_{(k)}$ is a subobject of $M\times M$, the synthetic 
language talks about $M_{(k)}$ {\em as if} it consisted of pairs of 
points of $M$; we shall for instance talk call such pair ``a pair of 
$k$'th order neighbours'' and write $x\sim _k y$ for $(x,y) \in 
M_{(2)}$.  For instance, the fact that $M_{(k)}$ is stable under the 
obvious twist map $M\times M \to M\times M$, we express by saying 
``$x\sim _k y$ implies $y\sim _k x$''.

 The ``set'' (scheme) of points $y\in M$ with $x\sim _k
y$, we also denote ${\cal M}_k (x)$, the {\em $k$'th order
neighbourhood},
or {\em $k$'th   monad}, around  $x$. The relation $\sim _k$ is reflexive and
symmetric, but not transitive; rather $x\sim _k y$ and $y\sim _l z$
implies $x\sim_{k+l}z$.  - Any map $f$ preserves these relations: 
$x\sim _k y$ implies $f(x)\sim _k f(y)$.

A quadratic differential form on $M$, i.e.  an element of $I^2 /I^3$,
can now be expressed: it is a function $g(x,y)$, defined whenever
$x\sim _2 y$, and so that $g(x,y)=0$ if $x\sim _1 y$.  If  further 
$g$ is positive definite, then we may {\em directly} think of 
$g(x,y)\in R$ as the {\em square distance} between $x$ and $y$.

For $M=R^n$, $M_{k}$ is canonically isomorphic to $M\times D_k (n)$:
$(x,y)\in R^n _{(k)}$ iff $y-x \in D_k (n)$; here, $D_k (n)$ is the
``infinitesimal'' scheme corresponding to a certain well known
Weil-algebra:

Recall that a Weil algebra is a finite dimensional ${\bf
R}$-algebra, where the nilpotent elements form a (maximal) ideal of 
codimension one.  The most basic Weil algebra is the ring of dual 
numbers
$${\bf R}[\epsilon ] = {\bf R}[X]/(X^2 ) = C^{\infty}({\bf R})/(x^2 );$$
the corresponding affine scheme is often denoted $D$, and is to be
thought of as a ``disembodied tangent vector''
(cf.\ Mumford \cite{mumford}, III.4, or Lawvere, \cite{fwl}).  The
reason is that maps of schemes $D\to M$ ($M$ a manifold, say) by
definition correspond to ${\bf R}$-algebra maps $C^{\infty}(M) \to
{\bf R}[\epsilon ]$, and such in turn correspond, as is known, to tangent
vectors of $M$.

Note that since ${\bf R}[\epsilon ] $ is a {\em quotient} algebra of
$C^{\infty}({\bf R})$, $D$ is, by the duality, a {\em sub}scheme of $R$;
this subscheme may be described synthetically as $\{ d\in R\mid d^2 =0 \}$, 
reflecting the fact that ${\bf R}[\epsilon ]$ comes about from 
$C^{\infty}({\bf R})$ by dividing out $x^2$.

More generally, for $k$ and $n$ positive integers, $D_k (n)$ is the
scheme corresponding to the Weil algebra which one gets from ${\bf R}[X_1 ,
\ldots ,X_n ]$ by dividing out by the ideal  generated by monomials of
degree $k+1$; or, equivalently, from $C^{\infty}({\bf R}^n )$ by the
ideal of functions that vanish to order $k+1$ at $\vo =(0,\ldots ,0)$ (it
is also known as the ``algebra of $k$-jets at $\vo$ in ${\bf R}^n$''). -- In
particular, $D_1 (1)$ is the ring of dual numbers described above.
Just as $D$ is the subscheme of $R$ described by $D= \{x\in R \mid
x^2 =0 \}$, $D_k (n)$ may be described in synthetic  language as 
$$\{(x_1 , \ldots ,x_n )\in R^n \mid x_{i_1}\cdot \ldots \cdot
x_{i_{k+1}}=0 \mbox{ for all } i_1 ,\ldots ,i_{k+1}\}.$$

The specific Weil algebras which form the algebraic backbone of the
present paper are the following (first studied for this purpose in
\cite{Laplace}).  For each
natural number $n\geq 2$, we consider the algebra $C^{\infty}({\bf
R}^n)/I_L$, where $I_L$ is the ideal generated by all $x_{i}^2 -x_{j}^2$
and all $x_i x_j$ where $i\neq j$.  The linear dimension of this algebra is
$n+2$; a basis may be taken to be (the classes mod $I_L$ of) the
functions $1, x_1 , \ldots ,x_n , x_1 ^2 + \ldots + x_n ^2$.  The
corresponding affine scheme we denote $D_L (n)$ or $D_L (R^n )$;  the
letter ``L'' stands for ``Laplace'', for reasons that will hopefully
become clear.  Using synthetic language, $D_L (n)$ may be described $$D_L (n) 
= \{(x_1 , \ldots , x_n )\in R^n \mid x_i ^2 = x_j ^2 ; x_i
x_j =0 \mbox{ for } i\neq j\}.$$

Note that $D_1 (n) \subseteq D_L (n) \subseteq D_2 (n)$.
The inclusion $D_1 (n) \subseteq D_L (n)$
corresponds to the quotient map $$C^{\infty}({\bf R}^n)/I_L
   \to C^{\infty}({\bf R}^n)/I_1$$
   which in turn comes about because $I_L \subseteq I_1$.  The kernel of
   this quotient map has linear dimension 1; a generator for it is the
   (class mod $I_L$ of) $x_1 ^2 + \ldots + x_n ^2$.

The following is  a {\em tautological} translation of  this fact:

\begin{prop}Any function $f:D_L (n) \to R$ which vanishes on $D_1 (n)$
is of the form $c\cdot ( x_1 ^2 + \ldots + x_n ^2 )$ for a unique $c\in R$.
\label{one}\end{prop}

\medskip The subscheme $D_k (n)\subseteq R^n$ can be described in
coordinate free terms; in fact, it is just the $k$-monad ${\cal M}_k
(\vo )$ around $\vo$.  More generally, for any finite dimensional vector 
space $V$, we can give an alternative description of ${\cal M}_k (\vo )$, 
which we also denote $D_k (V)$.  We only give this description for the 
case $k=1$ and $k=2$, which is all we need:

We have that $\vu \in D_1 (V)$ iff for any bilinear $B: V\times V \to
R$, $B(\vu ,\vu )=0$; this then also holds for any bilinear $V\times V
\to W$, with $W$ a finite dimensional vector space.  Similarly $\vu
\in D_2 (V)$ iff for any trilinear $C:V\times V\times V \to R$
, $C(\vu ,\vu ,\vu )=0$; this then also holds for any trilinear
$V\times V \times V \to W$, with $W$ a finite dimensional vector
space.

Any function $f: D_2 (V) \to W$ (with $W$ a finite dimensional
vector space) can uniquely be written in the form $\vu \mapsto f(\vo ) +
L(\vu ) + B(\vu ,\vu )$ with $L:V\to W$ linear and $B: V\times V \to
W$ bilinear symmetric.

If $V$ is equipped with a positive definite inner product, we shall
in the following Section also describe a subscheme $D_L (V)$ with 
$D_1 (V)\subseteq D_L (V)
\subseteq D_2 (V)$; for $V=R^n$ with standard inner product, it will 
be the $D_L (n)$ already described.

\section{L-neigbours in inner-product spaces}

For a 1-dimensional vector space $V$, we say that $a\in V$ is {\em
L-small} if it is 2-small, i.e. if $a\sim _2 0$.

Given an $n$-dimensional vector space $V$ ($n\geq 2$) with a positive
definite inner product $<-,->$.  We call a vector $\va \in V$ {\em 
L-small} if for all $\vu ,\vv \in V$ \begin{equation}<\va ,\vu ><\va ,\vv > = 
\frac{1}{n} <\va ,\va ><\vu ,\vv >.  \label{L-small}\end{equation}The ``set'' 
(scheme) of L-small vectors is denoted $D_L (V)$.  

It is clear that if $\va\in D_L (V)$, then  $\lambda \va \in D_L (V)$ for
any scalar $\lambda$. But $D_L (V)$ will not be stable under addition;
it is not hard to prove that if $\va$ and $\vb$ are L-small vectors, then
$\va + \vb$ is L-small precisely if
for all $\vu ,\vv \in V$ \begin{equation}<\va ,\vu ><\vb ,\vv > + <\va 
,\vv ><\vb ,\vu > = \frac{2}{n} <\va ,\vb ><\vu ,\vv >.
\end{equation}

\medskip

Let us analyze these notions for the case of $R^n$, with its standard
inner product. We claim

\begin{prop} The vector $\vt=(t_1 , \ldots ,t_n )$ belongs to $D_L
(R^n )$ iff
\begin{equation}t_1 ^2 = \ldots = t_n ^2\mbox{ ; and } t_i t_j =0
\mbox{ for } i\neq j .
\label{classic}\end{equation}
\end{prop}
(So $D_L (R^n )$ equals the $D_L (n)$ described above, or in \cite{Laplace} 
equation (8).)

\medskip

\noindent{\bf Proof.} If $\vt \in D_L (R^n )$, we have in particular for each 
$i=1, \ldots ,n$,
$$t_i ^2 = <\vt ,\ve _i ><\vt , \ve _i > = \frac{1}{n} <\vt ,\vt>,$$
where $\ve _1 , \ldots ,\ve _n $ is the standard (orthonormal) basis for
$R^n$. The right hand side here is independent of $i$. -- Also, if
$i\neq j$,
$$t_i t_j = <\vt ,\ve _i > <\vt , \ve _j > = \frac{1}{n} <\vt 
,\vt ><\ve _i , \ve _j > =0,$$
since $<\ve _i , \ve _j >=0$.

Conversely, assume that (\ref{classic}) holds.  Let $\vu$ and $\vv$ be 
arbitrary vectors, $\vu=(u_1 , \ldots ,u_n )$, and similarly for $\vv$. Then
$$<\vt ,\vu ><\vt ,\vv > = (\sum _i t_i u_i )(\sum _j t_j v_j )$$
$$=\sum_{i,j} t_i t_j u_i v_j = t_1 ^2 \sum _i u_i v_i ,$$
usin (\ref{classic}) for the last equality sign. But this is $t_1 ^2
<\vu ,\vv >$, and since, again by (\ref{classic}) $$t_1 ^2 = 
\frac{1}{n} (t_1 ^2 + \ldots + t_n ^2 ) = \frac{1}{n}
<\vt ,\vt >,$$
we conclude $<\vt ,\vu ><\vt ,\vv >= \frac{1}{n}<\vt ,\vt ><\vu ,\vv >$.

\medskip
As a Corollary, we get that for $\vv\in V$ (an $n$-dimensional
inner-product space), $\vv \in D_L (V) $ iff for some, or for any, 
orthonormal coordinate system for $V$, the coordinates of $\vv$ satisfy 
the equations (\ref{classic}). 

\medskip

From the coordinate characterization of $D_L (V)$ also immediately 
follows that $D_L (V) \subseteq D_2 (V)$.

\medskip
Here is an alternative characterization of L-small vectors, for inner
product spaces $V$ of dimension $\geq 2$ (the word ``self-adjoint''
may be omitted, but we shall need the Proposition in the form stated).
\begin{prop}The vector $\va$ belongs to $D_L (V)$ if and only if for
every self adjoint linear map $L:V\to V$ of trace zero, $<L(\va ), \va
>=0$
\label{linear-fugl}\end{prop} {\bf Proof.} We pick orthonormal coordinates,
and
utilize the ``coordinate'' description of $D_L (R^n )$.  Assume $\va
\in D_L (R^n )$, and assume $L$ is given by the symmetric matrix
$[c_{ij}]$ with $\sum c_{ii}=0$.  Then
$$<L(\va ),\va >= \sum _{ij} c_{ij}a_j a_i ;$$
since $a_i a_j =0$ if $i\neq j$, only the diagonal terms survive, and we get
$<L(\va ),\va > = \sum _i c_{ii} a_i ^2 = a_1 ^2 \sum c_{ii}$, since all 
the $a_i ^2$ are equal to $a_1 ^2$.  Since $\sum c_{ii} =0$, we get 0, 
as claimed.  Conversely, let us pick the $L$ given by the symmetric 
matrix with $c_{ij}=c_{ji}=1 (i\neq j)$, and all other entries $0$.  
Then
$$0=<L(\va ),\va > = a_i a_j + a_j a_i ,$$
whence $a_i a_j =0$.  Next let us pick the $L$ given by the matrix $c_{ii}=1$,
$c_{jj}=-1$ ($i\neq j$) and all other entries $0$.  Then
$$0=<L(\va ),\va > = a_i a_i - a_j a_j ,$$
whence $a_i ^2 = a_j ^2$.  So $\va \in D_L (R^n )$.

\medskip
We now consider the question of when a linear map $f:V\to W$ between
inner product spaces preserves L-smallness, i.e.\ when $f(D_L (V))
\subseteq D_L (W)$.  Let us call an $m\times n$ matrix {\em semi-conformal}
if the rows are mutually orthogonal, and have same (strictly positive)
square norm.  (This square norm is then called the {\em square
dilation} of the matrix, and is typically denoted $\Lambda$.) 
  The rank of a semi-conformal matrix is $m$, since its rows,
being orthogonal, are linearly independent.  It thus represents a
surjective linear map $R^n \to R^m$.

We have

\begin{prop} Let $f: V\to W$ be a surjective linear map between inner product
spaces.  Then t.f.a.e.\

1) $f(D_L (V))\subseteq D_L (W)$

2) In some, or any, pair of orthonormal bases for $V$, $W$, the matrix 
expression for $f$ is a semi-conformal matrix

In case these conditions hold, the common square norm $\Lambda$ of the 
rows of the matrix is characterized by: for all $\vz \in D_L (V)$
$$\frac{1}{m} <f(\vz ),f(\vz )> = \Lambda \frac{1}{n} <\vz ,\vz >,$$
(where $n=dim (V)$, $m=dim(W)$).
\label{linearconformal}\end{prop}
{\bf Proof.} Assume 1).  Pick orthonormal bases for $V$ and $W$,
thereby identifying $V$ and $W$ with $R^n$ and $R^m$, with standard inner
product.  Let the matrix for $f$ be $A= [ a_{ij}]$.  For all $\vz \in
D_L (V)$, we have by assumption that
$$(\sum _j a_{ij} z_j)^2 \mbox{ is independent of } i.$$
We calculate this expression:
\begin{equation}(\sum _j a_{ij} z_j)^2 = (\sum _j a_{ij} z_j) (\sum
_{j'} a_{ij'}
z_{j'}) =
\sum _j a_{ij}^2 z_j ^2 \label{a1}\end{equation} since the condition $\vz\in D_L
(V)$ implies that $z_j z_{j'} =0$ for $j\neq j'$, so all terms where $j\neq j'$ are
killed.  Also $z_j ^2 = z_1 ^2$, so bringing this factor outside the
sum, we get
\begin{equation} = z_1 ^2 (\sum _j a_{ij} ^2) = \frac{1}{n}
(\sum _k z_k ^2 ) (\sum _j a_{ij} ^2).\label{a2}\end{equation}
Since this is independent of $i$, then
so is $\sum _j a_{ij} ^2$, by the uniqueness assertion in Proposition
\ref{one}.  - The proof that the rows of $A$ are mutually orthogonal is 
similar (or see the proof for Theorem 3.2 in \cite{Laplace}).  -- 
Conversely assume 2), and assume that $\vz \in D_L (R^n )$.  We prove 
that $A\cdot \vz \in D_L (R^m )$.  The square of the $i$'th coordinate 
here is \begin{equation}(\sum _j a_{ij}z_j )^2 = z_1 ^2 \sum _j 
a_{ij}^2\label{a3}\end{equation} by the same calculation as before.  
But now the sum is independent of $i$, by assumption on the matrix 
$A$.  -- Similarly, if $i\neq i'$, the inner product of the $i$'th and 
$i'$'th row of $A\cdot \vz$ is $$(\sum _j a_{ij}z_j )(\sum _{j'} a_{i' 
j'}z_{j'}) =z_1 ^2
(\sum _j a_{ij}a_{i'j}),$$
using again the special equations that hold for the $z_j$'s; but now
the sum in the parenthesis is $0$ by the assumed orthogonality of the
rows of $A$.

Let $\Lambda$ be the common square norm of the rows of the matrix for
$f$.  Then for $\vz \in D_L (V)$, $$\frac{1}{m} <f(\vz ),f(\vz )> = \frac{1}{m}
\sum _i (\sum_j
a_{ij}z_j)(\sum _{j'} a_{ij'}z_{j'}),$$
and multiplying out, only the terms where $j=j'$ survive, since $\vz \in
D_L (V)$.  Thus we get
$$\frac{1}{m} \sum _i (\sum _j a_{ij}^2 z_j ^2 ) = \frac{1}{m} z_1 ^2
(\sum _i \sum _j a_{ij}^2)= \frac{1}{m} z_1 ^2
(\sum _i \Lambda ) $$
but this is $z_1 ^2 \Lambda$, since there are $m$ indices $i$.  On the
other hand, $z_1 ^2 = 1/n (\sum _j z_j ^2 )$.

\medskip

We have the following ``coordinate free'' version of Proposition
\ref{one} (derived from it by picking orthonormal coordinates):

\begin{prop}
Let $f_1 , f_2: D_L (V) \to R$ be functions which agree on $D_1 (V)$.
Then there exists a unique number $c\in R$ so that for all $\vz \in D_L
(M)$ we have
$$f_1 (\vz )- f_2 (\vz ) = c\cdot <\vz ,\vz >.$$
\label{little-lap}\end{prop}
Consider a map $f: D_2 (V) \to W$ with $f(\vo )=\vo$ and a symmetric 
bilinear $B: V\times V \to W$.  Let $b :V\to W$ denote the''quadratic 
`` map $\vu \mapsto B(\vu ,\vu )$.

\begin{lemma}The map $f$ takes $D_L (V)$ into $D_L (W)$ if and only if
$f + b$ does.
\label{little-b}\end{lemma}
{\bf Proof.} This is a simple exercise in degree calculus.  Assume $f$
has the property.  To prove that $f+b$ does, let $\va \in D_L
(V)$, and let $\vu ,\vv$ be arbitrary ``test'' vectors in $W$.
We consider $<f(\va ) + b(\va ) ,\vu ><f(\va )+b(\va ),\vv>$.  Using
bilinearity of inner product, this comes out as four terms, one of
which is $<f(\va ),\vu ><f(\va ),\vv )$, and three of which vanish for
degree reasons, thus for instance $<b(\va ),\vu ><f(\va ), \vv > =
<B(\va ,\va),\vu ><f(\va ), \vv >$ which contains $\va$ in a trilinear
way, so vanishes since $\va \in D_L (V) \subseteq D_2
(V)$.  So the left hand side in the test equation for L-smallness of 
$f(\va )+ b(\va )$ equals the left hand side in the test equation for 
L-smallness of $f(\va )$.  The right hand sides of the test equation is 
dealt with in a similar way.

\section{Riemannian metrics}
Recall from \cite{GCLCP}, \cite{Laplace} that a Riemannian metric $g$
on a manifold $M$ may be construed as an $R$-valued function defined
on the second neighbourhood $M_{(2)}$ of the diagonal, and vanishing
on $M_{(1)}\subseteq M_{(2)}$; we think of $g(x,y)$ as the {\em square
distance} between $x$ and $y$. Also $g$ should  be {\em positive
definite}, in a sense which is most easily expressed when passing to a 
coordinatized situation.  Since our arguments are all of completely 
local (in fact infinitesimal) nature, there is no harm
 in assuming that {\em one} chart covers all
of $M$, meaning that we have an embedding of $M$ as an open subset of
$R^n$, or of an abstract $n$-dimensional vector space $V$.  In this
case, each $T_x M$ gets canonically identified with $V$: to $\vu \in
V$, associate the tangent vector $t$ at $x$ given by $d\mapsto x +
d\cdot \vu$ for $d\in D$.  The vector $\vu$ is called the {\em
principal part} of $t$.  
  In this case $g$ is of the form
$$g(x,z)= G(x; z-x,z-x),$$
where $G:M\times V \times V \to R$ is bilinear symmetric in the two 
last arguments.  We require each $G(x;-,-)$ to be positive definite, 
i.e.\ $G(x;-,-)$ provides $V$ with an inner product (depending on 
$x$).  Since $T_x M$ is canonically identified with $V$, each $T_x M$ 
also acquires an inner product; this inner product can in fact be 
described in a coordinate free way, in terms of $g$ alone, cf.\ 
\cite{Laplace} formula (4).

\section{Symmetric affine connections, and the $\log$-$\exp$-bijection}
According to \cite{CTC}, an {\em affine connection} $\nabla $ on a 
manifold $M$ is a law $\nabla $ which allows one to complete any 
configuration (with $x\sim _1 y, x\sim _1 z$)

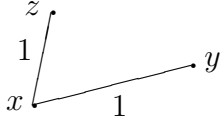
\vspace{5mm}

\begin{picture}(200,70)(-90,-10)
\put(20,6){\line(4,1){60}}
\put(20,6){\line(1,5){7}}
\put(10,4){$x$}
\put(85,21){$y$}
\put(17,40){$z$}
\put(50,2){1}
\put(14,23){1}
\put(19,5){{\bf .}}
\put(79,20){{\bf .}}
\put(26,40){{\bf .}}
\end{picture}

\noindent into a configuration

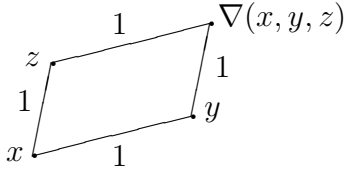
\vspace{5mm}

\begin{picture}(200,70)(-90,-10)
\put(20,6){\line(4,1){60}}
\put(20,6){\line(1,5){7}}
\put(80,21){\line(1,5){7}}
\put(27,41){\line(4,1){60}}
\put(10,4){$x$}
\put(85,21){$y$}
\put(17,40){$z$}
\put(90,55){$\nabla(x,y,z)$} \put(19,5){{\bf .}} \put(79,20){{\bf .}}
\put(26,40){{\bf .}} \put(86,55){{\bf .}}

\put(50,2){1}
\put(50,52){1}
\put(14,23){1}
\put(89,36){1}

\end{picture}

\vspace{5mm}

\noindent (with $z\sim _1 \nabla (x,y,z)\sim _1 y$), to be thought of
as an ``infinitesimal parallelogram according to $\nabla$''.  There is only
one axiom assumed:
$$\nabla (x,x,z)=z;\; \nabla (x,y,x)=y.$$
     If $\nabla (x,y,z)=\nabla (x,z,y)$ for all $x\sim _1 y$, $x\sim _1 z$, we
     call the connection {\em symmetric}.

     In a coordinatized situation, i.e. with $M$ identified with an open
     subset of a finite dimensional vector space $V$, the data of an
     affine connection $\nabla$ may be encoded by a map $\Gamma : M\times
     V \times V \to V$, bilinear in the two last arguments, namely
     $$\nabla (x,y,z)= y-x+z + \Gamma (x; y-x,z-x),$$ so that $\Gamma $
     measures the discrepancy between ``infinitesimal parallelogram
     formation according to $\nabla$'' and the corresponding
     parallelograms according to the affine
     structure of the vector space $V$. This $\Gamma$ is
     the ``union of'' the Christoffel symbols; and $\nabla$ is symmetric iff 
     $\Gamma (x; -,-)$ is.

A fundamental result in differential geometry is the existence of the
Levi-Civita connection associated to a Riemann metric $g$. This
result can be formulated synthetically, without reference to tangent bundles or 
coordinates, namely: given a Riemann metric $g$ on a manifold, then 
there exists a unique symmetric connection $\nabla$ on $M$ with the 
property that for any $x\sim _1 y$, the map $\nabla (x,y,-):{\cal M}_1 
(x)\to {\cal M}_1 (y)$ preserves $g$, i.e.  for $z\sim _1 x$, $u\sim 
_1 x$,
$$g(\nabla (x,y,z),\nabla (x,y,u))= g(z,u).$$

(This latter condition is equivalent to: the differential of $\nabla
(x,y,-)$ at $x$ is an inner-product preserving linear map $T_x M \to T_y M$.)

\medskip There is, according to \cite{KL} Theorem 4.2, an alternative way of 
encoding the data of a symmetric affine connection on $M$, namely as a 
``partial exponential map'', meaning a bijection (for each $x\in M$) 
${\cal M}_2 (x)\cong D_2 (T_x M)\subseteq T_x (M)$, with certain 
properties.  We describe how such bijection $\exp _x :D_2 (T_x M) \to 
{\cal M}_2 (x)$ is related to the connection $\nabla$ (and this 
equation characterizes $\exp _x$ completely):
$$\exp _x ((d_1 + d_2)t)= \nabla (x, t(d_1 ), t(d_2 )),$$
where $t\in T_x M$ and $d_1 , d_2 \in D$ (this implies $(d_1 + d_2 )t
\in D_2 (T_x M)$).

Since $\nabla (x,y,x)=y$, it follows by taking $d_2 =0$ that $\exp(
d_1 t) = t(d_1 )$, so that the partial exponential map ${\cal M}_2 (0)
\to {\cal M}_2 (x) $ is an extension of the ``first order'' partial
exponential map ${\cal M}_1 (0)
\to {\cal M}_1 (x) $, as considered in \cite{Krynica}; the first
order exponential map is ``absolute'' in the sense that its
construction does not depend on a metric $g$ on $M$.

In the coordinatized situation with $M\subseteq V$ an open subset of
a vector space $V$,
     the second order exponential map corresponding to $\nabla$ is given as
      follows.  Note first that since now $M$ is an open subset of $V$, $T_x
      (M)$ may be identified with $V$ canonically, via the usual notion of
      ``principal part'' of a tangent vector to $V$.  Let $\vu \in D_2 (V)$. Then
      $$\exp _x (\vu ) = x +\vu + \frac{1}{2} \Gamma (x; \vu, \vu ).$$
      This
      is an element in $M \subseteq V$, since $M$ is open, in fact, it is an
      element of ${\cal M}_2 (x)$.

The inverse of $\exp _x$ we of course have to call $\log _x$; in the
coordinatized situation $M\subseteq V$, it is given as follows:
let $y\sim
_2 x$; then $y= x+\vu$ with $\vu \in D_2 (V)$, and
     $$\log _x (x+\vu ) = \vu - \frac{1}{2}\Gamma (x; \vu , \vu ).$$
The fact that the map $\log _x$ thus described is inverse for $\exp _x$ is a
simple calculation using bilinearity of $\Gamma (x; -,-)$, together
with $\Gamma (x; \vu , \Gamma (x; \vu ,\vu )) =0$, and $\Gamma (x;
\Gamma (x;\vu ,\vu ) , \Gamma (x; \vu ,\vu ))=0$, and these follow
because they are trilinear (respectively quatrolinear) in the
arguments where $\vu$ is substituted.

--The following gives an ``isometry'' property of the 
$\log$-$\exp$-bijection.  (It does not depend on the relationship 
between the metric $g$ and the affine connection/partial exponential.)

\begin{prop} For $z\sim _2 x$, $g(x,z)= <\log _x z, \log _x z>$.
\label{exp-isometry}\end{prop}
{\bf Proof.} We work in a coordinatized situation $M\subseteq V$, so 
that $g$ is encoded by $G: M\times V \times V \to R$, and the connection 
is encoded by $\Gamma :M\times V \times V \to V$, with both $G$ and 
$\Gamma$ bilinear in the two last arguments.  Let $z\sim _2 x$, so $z$ 
is of the form $x+\vu$ with $\vu \in D_2 (V)$.  Then on the one hand
$$g(x,z)= G(x;\vu ,\vu ),$$
and on the other hand,
$\log _x (z) = \vu -1/2 \Gamma (x; \vu ,\vu )$ so that
$$<\log _x z, \log _x z > = G(x; \vu -\frac{1}{2}\Gamma (x;\vu ,\vu ),
\vu -\frac{1}{2}\Gamma (x;\vu ,\vu )),$$
and expanding this by bilinearity, we get $G(x; \vu ,\vu )$ plus some 
terms which vanish because they are tri- or quatro-linear in $\vu$.

\section{Mirror image} Using the (second order) partial exponential
map, we can give a simple description of the infinitesimal symmetry
(\cite{Laplace}) which any Riemannian manifold has.  Let $z\sim _2 x$
in $M$.  Its {\em mirror image } $z'$ in $x$ is defined by
$$z' := \exp_ x (- \log _x (z)).$$

In the coordinatized situation $M\subseteq V$, we can utilize the
formulae for $\log$ and $\exp$ given in terms of $\Gamma$ to get the
following formula for mirror image formation. If $z = x + \vu$ with
$\vu \in D_2 (V)$, we get
$$z' = x - \vu + \Gamma (x;\vu , \vu ).$$
This is a calculation much similar to the one above, namely, cancelling
terms of the form $\Gamma (x; \Gamma (x;\vu , \vu), \vu )$ or $\Gamma
(x; \Gamma (x;\vu , \vu), \Gamma (x;\vu , \vu) )$, these being tri- or
quatro-linear in $\vu$.  A similar calculation will establish that
$z'' = z$.

Note also that if $\vu \in D_1 (V)$, and $z=x+\vu$, then $z' = x 
-\vu$.

    From this follows

\begin{lemma} Given $x\in M$. Let $f:M\to R$. The function
$\tilde{f}:{\cal M}_2 (x)\to R$ defined by
$$\tilde{f} (z)= f(z' ) + f(z) -2f(x)$$
vanishes on ${\cal M}_1 (x)$.
\label{van}\end{lemma} For, if $df$ denotes the differential of $f$ at $x$,
and $z=x+\vu$ with $\vu \in D_1 (V)$, the right hand side here
is
$$(f(x) +df(-\vu )) +( f(x)+df (\vu)) -2f(x),$$
and this is $0$ since $df$ is linear.

\section{L-neighbours in a Riemannian manifold} We consider a
Riemannian manifold $(M, g)$, and the various structures on $M$ derived
from it, as in the previous sections. In particular, we have the
partial exponential map $\exp$, and its inverse $\log$.  Using these maps, 
we shall transport the L-neighbour relation from the inner-product 
spaces $T_x M$ back to a relation in $M$.  Explicitly, \begin{defn}Let 
$x\sim _2 z$ in $M$.  We say that $x\sim _L z$ if $\log _x (z)$ is 
L-small in the inner product space $T_x M$ (with inner product derived 
from $g$).  \end{defn}Note that this is not apriori a symmetric 
relation, since $\log (x,z)$ and $\log (z,x)$ are not immediately 
related -- they belong to two different vector spaces $T_x M $ and 
$T_z M$; in a coordinatized situation $M\subseteq V$, both these 
vector spaces may be canonically identified with $V$, but the notion 
of $\exp$ and $\log$ depend on inner products, and $V$ in general gets 
different inner products from $T_x M$ and $T_z M$.  In \cite{Laplace}, 
the question of symmetry of the relation $\sim _L$ was left open (and 
the relation $\sim _L$ was defined in a different, more complicated way).  We 
state without proof:

\begin{prop}The L-neighbour relation is symmetric.
\end{prop}
This fact will not be used in the present paper. It depends on the fact
that parallel transport according to $\nabla$ preserves L-smallness,
being an isometry.

\medskip
The following is the fundamental property of L-neighbours, and
provides the link to the Laplace operator and harmonic functions, and
more generally to harmonic morphisms. It is identical to Theorem 2.4 in 
\cite{Laplace}, but the argument we give presently is more canonical 
(does not depend on chosing a geodesic coordinate system):

\begin{thm} For any $f:{\cal M}_L (x) \to R$, there exists a unique
number $c$ so that for all $z\in {\cal M}_L (x)$,
\begin{equation}f(z)+f(z')-2f(x)= c \cdot
g(x,z).\label{fund}\end{equation}
\label{big-lap}\end{thm}
{\bf Proof.} Consider the composite of $\exp _x$ with the function 
$\tilde{f}$ of $z$ described by the left hand side of (\ref{fund}),
$$\begin{diagram}
D_L (T_x M)& \rTo^{\exp _x} & {\cal M}_L (x) & \rTo^{\tilde{f}} & R
\end{diagram}$$
It is a function defined on $D_L (T_x M ) \subseteq T_x M$.  It follows 
from Lemma \ref{van}  that this function vansihes on $D_1 (T_x M )$, and 
thus is constant multiple of the square-norm function $T_x M \to R$, 
by Proposition \ref{little-lap},
$$\tilde{f}(\exp _x (\vu ,\vu )) = c\cdot <\vu ,\vu >.$$
Apply this to $\vu = \log _x z$ for $z\sim _L x$;
we get
$$\tilde{f}(z) = \tilde{f}(\exp _x (\log _x (z))) = c\cdot <\log _x z 
,\log _x z>,$$
which is $c\cdot g(x,z)$ by Proposition \ref{exp-isometry}.

\medskip
 For any  function $f:M \to R$, we can for each $x\in M$ consider the 
corresponding $c$, characterized by (\ref{fund}); this gives a 
function $c: M\to R$, and we define $\Delta (f)$ to be $n$ times this 
function, in other words, the function $\Delta (f)$ is characterized 
by: for each pair $x \sim _L z$
\begin{equation}f(z)+f(z')-2f(x) = \frac{\Delta (f)(x)}{n}
g(x,z),\label{Delta-def}\end{equation}
where $z'$ denotes the mirror image of $z$ in $x$.
(This $\Delta$ operator can be proved to be the standard Laplace
operator, cf.
\cite{Laplace}.)

We call $f$ a {\em harmonic} function if $\Delta (f) =0$. Thus
harmonic functions are characterized by the {\em average value
property}:
for any $x\sim _L z$, $f(x)$ is the average of $f(z)$ and $f(z')$.
This property can also be expressed: for any $z\sim _L x$, $f(z')$ is
the mirror image of $f(z)$ in $f(x)$, where mirror image of $b$ in
$a$ for $a,b \in R$ means $2a -b$. This is also the mirror-image
formation in $R$ w.r.to the standard Riemannian metric given by
$g(a,b)=(b-a)^2$.

This observation prompts the following definition:

\begin{defn}Let $(M,g)$ and $(N,h)$ be Riemannian 
manifolds, and let
$\phi: M\to N$ a map. We say that $\phi$ is a {\em harmonic map} if
it preserves mirror image formation of L-neighbours $x,z$,
$$\phi (z')=\phi (z)',$$
where the prime denotes mirror image formation in $x$ w.r.to $g$ and
in $\phi (x)$ w.r.to $h$, respectively.
\end{defn}
Note that even if $z$ is an L-neighbour of $x$,
$\phi (z)$ may not be an L-neighbour of $\phi (x)$, but it will
  be a 2-neighbour of $\phi (x)$, so that the notion of mirror image of it 
  makes sense.  -- The notion may be localized at $x$: $\phi$ is a 
  harmonic map {\em at} $x$ if for all $z\sim _L x$, $\phi (z')= 
  \phi(z)'$.)

A stronger notion than harmonic {\em map} is that of harmonic {\em morphism};
this is a map
which is as well a harmonic map, and  is also {\em semi- (or horizontally)
conformal} in the sense of the next section. (The terminology is not
very fortunate, but
classical, cf.\  \cite{BW}.)

\section{Semi-conformal maps} We consider again two Riemannian manifolds
$(M,g)$ and $(N,h)$, and a submersion $\phi :M\to N$. It defines
a ''vertical'' foliation, whose leaves are the (components of) the
fibres of $\phi$, and hence the transversal distribution consisting
of $Ker (df_x )^{\perp} \subseteq T_x M$.
(This ``horizontal''
distribution can also be described in purely combinatorial terms 
without reference to the tangent bundle.)

Recall (from \cite{BW}, say) that $\phi$ is called {\em
semi-conformal} (or horizontally conformal) at $x\in M$,
with square-dilation $\Lambda >0$, if the
linear map
$df_x : T_x M \to T_{\phi (x)}N$ is semi-conformal
with square-dilation $\Lambda >0$, in the sense of Section 2.
(This property can also be expressed combinatorially.)
 The following is 
a generalization of Theorem 3.2 in \cite{Laplace} (which dealt with 
the case of a diffeomorphism $\phi$).  \begin{thm}Let $\phi :M\to N$ 
be a submersion, and let $x\in M$.  Then t.f.a.e.:

1) $\phi$ is semi-conformal at $x$ (for some $\Lambda >0$)

2) $\phi$ maps ${\cal M}_L (x)$ into ${\cal M}_L (\phi (x))$.
\label{semi-c}\end{thm}
{\bf Proof.} Consider the diagram
\begin{diagram}
{\cal M}_2 (x)& \rTo ^{\phi} & {\cal M}_2 (\phi (x))\\
\dTo ^{\log _x}&&\dTo _{\log _{\phi (x)}}\\
D_2 (T_x M)& \pile{\rTo^f \\\rTo _{ d\phi _x}} & D_2 (T_{\phi
(x)}
N)
\end{diagram}
where $f$ is the unique map making the diagram commutative, and where
$d\phi _x$ is (the restriction of) the differential of $\phi$.  It
does not make the diagram commutative, but when restricted to ${\cal
M}_1 (x)$, it does, by the very definition of differentials.  So $f$
and $d\phi _x$ agree on ${\cal
M}_1 (x)$, and hence differ by a quadratic map $b$.  It then follows
from Lemma \ref{little-b}
that $f$ maps $D_L (T_x M)$ into $D_L (T_{\phi (x)})
$if and only if $d\phi _x$ does.  By definition, ${\cal
M}_L
(x)$ comes about from $D_L (T_x M)$ by transport along the
$\log$-$\exp$-bijection, so $\phi$ preserves ${\cal M}_L$ iff $f$
preserves $D_L$.
On the other hand, by the Proposition \ref{linearconformal}, semi-conformality
of $d\phi _x$ is equivalent to $d\phi _x$ preserving $D_L$.

\medskip
We may summarize the results of the last two sections by stating the
following (which may be taken as definitions of these notions, but
couched in purely geometric/combinatorial language): let $\phi :M \to
N$ be a submersion between Riemannian manifolds. Then

\begin{itemize}
\item $\phi$ is a {\em harmonic map} if it preserves mirror image
formation of L-neighbours

\item $\phi$ is a {\em semi-conformal map} if it preserves the notion
of L-neighbour

\item $\phi$ is a {\em harmonic morphism} if it has both these
properties.
\end{itemize}

If the codomain is $R$, any 2-neigbour is an L-neighbour, so any map
to $R$ is automatically semi-conformal, so for codomain $R$, harmonic
map and harmonic morphism means the same thing.  Such a map/morphism is 
in fact exactly a harmonic function $M\to R$.

All three notions make sense ``pointwise'': $\phi$ is a harmonic {\em at}
$x\in M$ if it preserves mirror image formation of L-neighbours of
$x$. For this to make sense, we don't need $\phi$ to be defined on all
of $M$, because the property only depends on the 2-jet of $\phi$ at 
$x$, meaning the restriction of $\phi$ to ${\cal M}_2 (x)$.

\section{Sufficiency of harmonic 2-jets}

By 2-jets, we understand in this Section  2-jets of $R$-valued functions; so 
a 2-jet at $x\in M$ is a map ${\cal M}_2 (x) \to R$.  If $M$ is a 
Riemannian  manifold, we say that such a 2-jet $f$ 
is {\em harmonic} if it preserves mirror image formation of 
L-neigbours of $x$, $f(z')= 2f(x) - f(z)$, for all $z\sim _L x$.

Among such harmonic 2-jets, we have in particular those of the form
\begin{equation}
\begin{diagram}{\cal M}_2 (x) &\rTo^{\log _x}& T_x M &\rTo^{p}& R,
\end{diagram}\label{affine}\end{equation}
where the last map $p$ is linear.  For, by construction of mirror
image in terms of $\log _x$,
$\log _x (z')= -\log _x (z)$, and this mirror image formation is
preserved by $p$ (here, we don't even need $z\sim _L x$, just $z\sim _2 x$).

Another type of harmonic 2-jet are those of the form
\begin{equation}
\begin{diagram}
{\cal M}_2 (x) &\rTo^{\log _x}& T_x M &\rTo^q & R
\end{diagram}\label{tracezero}\end{equation}
where $q$ is a ``quadratic map of trace 0'', meaning
$q(\vu ) = <L(\vu ),\vu >$ for some selfadjoint $L:T_x M \to T_x M$ of
trace zero.  For, $z\sim _L x$ means by definition that $\log _x (z)
\in D_L (T_x M)$, and quadratic trace zero maps kill $D_L$, by 
Proposition \ref{linear-fugl}.

These two special kinds of harmonic jets are the only ones that we
shall use in the proof of the following ``recognition Lemma'':

\begin{lemma}There are sufficiently many harmonic 2-jets to recognize
mirror image formation in $x$, and to recognize L-neighbours of $x$.
\label{recognition}\end{lemma}
Precisely, if $z$ and $\tilde{z}$ are 2-neighbours of $x$, and
$f(\tilde{z})= 2f(x) - f(z)$ for all harmonic 2-jets $f$, then
$\tilde{z}=z'$; and if $z$ is a 2-neighbour of $x$ such that $f(z)=0$
for all harmonic 2-jets $f$ which vanish on ${\cal M}_1 (x)$, then
$z\sim _L x$.

\medskip
\noindent {\bf Proof.} The first assertion follows because $\log _x 
(z') = -\log _x (z)$, and because there are sufficiently many linear 
$p:T_x M\to R$ to distinguish any pair of vectors ($T_x M$ being 
finite-dimensional).  The second assertion follows because $\log _x$ 
maps ${\cal M}_L (x)$ bijectively onto $D_L (T_x M)$, and the latter 
is recognized by quadratic trace zero maps, by Proposition 
\ref{linear-fugl}.

\medskip
There is a partial converse:

\begin{prop}
Let $f: {\cal M}_2 (x)\to R$ be a harmonic 2-jet which vansihes at 
${\cal M}_1 (x)$.  Then it vanishes at ${\cal M}_L (x)$.\label{converse}\end{prop}
{\bf Proof.} Let $b$ denote the composite $f\circ \exp _x :D_2 (T_x M) 
\to R$.  The vanishing assumption on $f$ implies that there is a unique 
quadratic map $T_x M\to R$ extending $b$.  It suffices to prove that 
$b(\vu )=0$ for any $\vu \in D_L (T_x M)$.  Let $z$ denote $\exp _x (\vu 
)$; then $z\in {\cal M}_L (x)$.  Harmonicity of $f$ at $x$ implies 
$f(z)+f(z') =0$ by Theorem \ref{big-lap}, and hence $b(\vu )+ b(-\vu 
)=0$.  But $b$ is a even function, being quadratic, hence $b(\vu )=0$.

\section{Characterization Theorem}
The following Theorem is now  almost immediate in view of the
combinatorial/geometric description of harmonic maps and  semi-conformal
maps.  It is a version of the Characterization Theorem of Fuglede and 
Ishihara, cf.\ \cite{BW} Theorem 4.2.2.

\begin{thm}
Given a submersion $\phi :M \to N$ between Riemannian manifolds, and
let $x\in M$. Then t.f.a.e.\

1) $\phi$ is a harmonic morphism at $x$

2) for any harmonic 2-jet $f$ at $\phi (x)$, $f\circ \phi : M\to R$
is a harmonic 2-jet.

\label{three}\end{thm}
(The Theorem in the classical form talks about harmonic {\em germs}
at $\phi (x)$, rather than harmonic 2-{\em jets}. The ``upgrading'' of our
version to the classical one thus depends on a rather deep existence
theorem: any harmonic 2-jet comes about by restriction from a
harmonic germ, see Appendix of \cite{BW}. Such existence results are
beyond the scope of our methods.)

\medskip \noindent {\bf Proof.} Assume that $\phi$ is a harmonic morphism at 
$x$, and let $f$ be a harmonic 2-jet.  Let $z\sim _L x$.  Then $\phi 
(z')=(\phi (z))'$, since $\phi$ is a harmonic map; also $\phi (z) \sim 
_L \phi (x)$ since $\phi$ is semi-conformal.  So $f$ preserves the 
mirror image of $\phi (z)$.  So both $\phi$ and $f$ preserve the 
relevant mirror images, hence so does the composite $f\circ \phi : 
{\cal M}_2 (x) \to R$; hence it is a harmonic 2-jet.

Conversely, suppose $\phi$ has $f\circ \phi$ harmonic for all harmonic
2-jets $f$ at $\phi (x)$. Let $z\sim _L x$. To prove
$\phi(z')=(\phi(z))'$, it suffices, by the Recognition Lemma (applied
to $N$) to prove that all harmonic 2-jets $f$ at $\phi (x)$ take $\phi
(z')$ to the mirror image of $\phi (z)$. But by assumption $f\circ
\phi$ is harmonic at $x$, so preserves mirror image. -- Also, to
prove $\phi (z)\sim _L \phi (x)$, it suffices by the Recognition
Lemma to prove that any harmonic 2-jet at $\phi (x)$, vanishing on
${\cal M}_1 (\phi (x))$, kills $\phi (z)$. But by assumption, $f\circ
\phi$ is a harmonic 2-jet, and it vanishes at ${\cal M}_1 (x)$, so by 
Proposition \ref{converse}, it kills $z$.  So 
$f(\phi (z))=0$, so $\phi (z)\sim _L \phi (x)$.  This proves the 
Theorem.

\end{document}